\documentclass[12pt,a4paper]{article}
\usepackage[utf8]{inputenc}
\usepackage{mathrsfs,amsthm,xcolor,verbatim,bbm,amsmath,amsfonts,amssymb,nicefrac,enumitem,hyperref,bm,mathtools,xparse,etoolbox,graphicx,subcaption}
\usepackage[a4paper, left=2cm, right=2cm]{geometry}
\usepackage{algorithm,algpseudocode}
\usepackage[nocompress]{cite}

\usepackage[capitalise,sort]{cleveref}

\crefname{enumi}{item}{items}
\crefname{equation}{}{}
\crefname{figure}{Figure}{Figures}
\crefname{listing}{Source code}{Source codes}
\crefname{lstlisting}{Source code}{Source codes}
\crefname{cor}{Corollary}{Corollaries}

\crefname{subsection}{Subsection}{Subsections}

\hypersetup{
    colorlinks,
    linkcolor={blue!80!black},
    citecolor={green},
    urlcolor={blue!80!black}
}

\usepackage{tcolorbox}
\tcbuselibrary{skins, breakable, listings}  
\newcounter{algorithmcounter}  
\renewcommand{\thealgorithmcounter}{\arabic{algorithmcounter}}  
\newtcolorbox[use counter=algorithmcounter]{myalgorithm}[3][]{
	enhanced,
	breakable,
	fonttitle=\bfseries,
	title=Algorithm~\thealgorithmcounter: #2,
	label={#3},  
	label type=algorithmcounter,
	#1,
	colframe=black,
	colback=white,
	coltitle=black,
	colbacktitle=white,
	sharp corners,
	boxrule=0.5pt,
	boxsep=1mm,
	top=1mm,
	bottom=1mm,
	left=0mm,
	right=0mm
}

\newcommand{\mycomment}[1]{\hfill \textcolor{gray}{\textit{\# #1}}}

\crefname{algorithmcounter}{Algorithm}{Algorithms}
\Crefname{algorithmcounter}{Algorithm}{Algorithms}
\crefname{line}{line}{lines}
\Crefname{line}{Line}{Lines}

\newcommand{\defaultParamDim}{d}
\newcommand{\momentprocess}{{\bf m}}
\newcommand{\Momentprocess}{{\bf v}}
\newcommand{\numberofsteps}{N}
\newcommand{\timesteps}{M}
\newcommand{\nos}{n}
\newcommand{\AvgPx}{\chi}
\newcommand{\AvgParam}{A}
\newcommand{\AvgPb}{K}
\newcommand{\AvgPc}{\mathbb{A}}
\newcommand{\avgPa}{\delta}

\newcommand{\defaultdataDim}{\mathscr{d}}

\newcommand{\defaultStochLoss}{\mathscr{l}}
\newcommand{\defaultStochGradient}{\mathscr{g}}

\hypersetup{
	colorlinks,
	linkcolor={blue!80!black},
	citecolor={green},
	urlcolor={blue!80!black}
}

\theoremstyle{plain}
\newtheorem{theorem}{Theorem}[section]

\newtheorem{prop}[theorem]{Proposition}

\newtheorem{definition}[theorem]{Definition}

\DeclareMathAlphabet{\mathpzc}{OT1}{pzc}{m}{it}

\DeclareFontEncoding{LS1}{}{}
\DeclareFontSubstitution{LS1}{stix}{m}{n}
\DeclareMathAlphabet{\mathscr}{LS1}{stixscr}{m}{n}

\newcommand{\E}{\mathbb{E}}
\renewcommand{\P}{\mathbb{P}}

\newcommand{\R}{\mathbb{R}}
\newcommand{\N}{\mathbb{N}}

\newcommand{\ssum}{\textstyle\sum}

\DeclarePairedDelimiter{\norm}{\lVert}{\rVert}
\DeclarePairedDelimiter{\abs}{\lvert}{\rvert}
\DeclarePairedDelimiter{\rbr}{(}{)}
\DeclarePairedDelimiter{\br}{[}{]}
\DeclarePairedDelimiter{\cu}{\{}{\}}
\DeclarePairedDelimiter{\spro}{\langle}{\rangle}

\renewcommand{\d}{ \mathrm{d}}

\newcommand{\andq}{\text{and}\qquad}

\newcommand{\realization}[1]{\mathscr{N}^{ #1 }}

\ExplSyntaxOn

\bool_new:N \g_noteobserve

\NewDocumentCommand{\setnote}{}{
	\bool_gset_true:N \g_noteobserve
}

\NewDocumentCommand{\setobserve}{}{
	\bool_gset_false:N \g_noteobserve
}

\NewDocumentCommand{\nobs}{ o }{
	\IfValueT{#1}{
		\str_if_eq:noTF {note} {#1} {
			\bool_gset_true:N \g_noteobserve
		} {
			\str_if_eq:noTF {Note} {#1} {
				\bool_gset_true:N \g_noteobserve
			} {
				\bool_gset_false:N \g_noteobserve
			}
		}
	}
	\bool_if:nTF { \g_noteobserve } {
		\bool_gset_false:N \g_noteobserve
		note
	} {
		\bool_gset_true:N \g_noteobserve
		observe
	}
	\IfValueF{#1}{~}
}

\NewDocumentCommand{\Nobs}{ o }{
	\IfValueT{#1}{
		\str_if_eq:noTF {note} {#1} {
			\bool_gset_true:N \g_noteobserve
		} {
			\str_if_eq:noTF {Note} {#1} {
				\bool_gset_true:N \g_noteobserve
			} {
				\bool_gset_false:N \g_noteobserve
			}
		}
	}
	\bool_if:nTF { \g_noteobserve } {
		\bool_gset_false:N \g_noteobserve
		Note
	} {
		\bool_gset_true:N \g_noteobserve
		Observe
	}
	\IfValueF{#1}{~}
}

\int_new:N \g_furthermore

\NewDocumentCommand{\Moreover}{ o o }{
	\IfValueT{#1}{
		\str_case:nn {#1} {
			{Furthermore} {\int_set:Nn {\g_furthermore} {0}}
			{Moreover} {\int_set:Nn {\g_furthermore} {1}}
			{In~addition} {\int_set:Nn {\g_furthermore} {2}}
			{note} {\bool_gset_true:N \g_noteobserve}
			{observe} {\bool_gset_false:N \g_noteobserve}
		}
		\IfValueT{#2}{
			\str_case:nn {#2} {
				{Furthermore} {\int_set:Nn {\g_furthermore} {0}}
				{Moreover} {\int_set:Nn {\g_furthermore} {1}}
				{In~addition} {\int_set:Nn {\g_furthermore} {2}}
				{note} {\bool_gset_true:N \g_noteobserve}
				{observe} {\bool_gset_false:N \g_noteobserve}
			}
		}
	}
	\int_case:nn { \int_mod:nn {\g_furthermore} {3} } {
		{ 0 } { Furthermore,~\nobs that}
		{ 1 } { Moreover,~\nobs that}
		{ 2 } { In~addition,~\nobs that}
	}
	\int_incr:N \g_furthermore
	\IfValueF{#1}{~}
}

\bool_new:N \g_hencetherefore

\NewDocumentCommand{\hence}{}{
	\bool_if:nTF { \g_hencetherefore } {
		\bool_gset_false:N \g_hencetherefore
		hence~
	} {
		\bool_gset_true:N \g_hencetherefore
		therefore~
	}
}

\NewDocumentCommand{\Hence}{}{
	\bool_if:nTF { \g_hencetherefore } {
		\bool_gset_false:N \g_hencetherefore
		Hence,~we~obtain~
	} {
		\bool_gset_true:N \g_hencetherefore
		Therefore,~we~obtain~
	}
}

\seq_new:N \g_cflist_loaded
\seq_new:N \g_cflist_pending

\NewDocumentCommand{\cfadd}{ m }
{
	\seq_if_in:NnF \g_cflist_loaded { #1 } {
		\seq_if_in:NnF \g_cflist_pending { #1 } {
			\seq_gput_right:Nn \g_cflist_pending { #1 }
		}
	}
}

\NewDocumentCommand{\cfconsiderloaded}{ m }{
	\seq_gput_right:Nn \g_cflist_loaded {#1}
}

\NewDocumentCommand{\cfremove}{ m }
{
	\seq_gremove_all:Nn \g_cflist_pending { #1 }
}

\NewDocumentCommand{\cfload}{ o }
{
	\seq_if_empty:NTF \g_cflist_pending {\unskip} {
		(cf.\ \cref{\seq_use:Nn \g_cflist_pending {,}})\IfValueTF{#1}{#1~}{\unskip}
		\seq_gconcat:NNN \g_cflist_loaded \g_cflist_loaded \g_cflist_pending
		\seq_gclear:N \g_cflist_pending
	}
}

\NewDocumentCommand{\cfclear} {} {
	\seq_gclear:N \g_cflist_loaded
	\seq_gclear:N \g_cflist_pending
}

\NewDocumentCommand{\cfout}{ o }
{
	\seq_if_empty:NTF \g_cflist_pending {\unskip} {
		(cf.\ \cref{\seq_use:Nn \g_cflist_pending {,}})\IfValueTF{#1}{#1~}{\unskip}
		\seq_gclear:N \g_cflist_pending
	}
}

\NewDocumentCommand{\ifnocf} { m } {
	\seq_if_empty:NT \g_cflist_pending { #1 }
}

\ExplSyntaxOff

\NewDocumentEnvironment{cproof}{m}
{\begin{proof}[Proof of \cref{#1}]}%
	{\noindent The proof of \cref{#1} is thus complete.
\end{proof}}

\NewDocumentEnvironment{cproof2}{m}
{\begin{proof}[Proof of \cref{#1}]}%
	{\noindent This completes the proof of \cref{#1}.
\end{proof}}

\hypersetup{colorlinks=true}

\makeatletter
\@namedef{subjclassname@2020}{%
	\textup{2020} Mathematics Subject Classification}
\makeatother

\ExplSyntaxOn

\bool_new:N \g_forexample

\NewDocumentCommand{\eg}{ o }{
\IfValueT{#1}{
\str_if_eq:noTF {fe} {#1} {
\bool_gset_true:N \g_forexample
} {\bool_gset_false:N \g_forexample}
}
\bool_if:nTF { \g_forexample } {
\bool_gset_false:N \g_forexample
for~example
}{
\bool_gset_true:N \g_forexample
for~instance
}
}

\ExplSyntaxOff

\makeatletter
\ExplSyntaxOn
\seq_new:N \g_abbrs
\prop_new:N \g_abbr_counts
\tl_new:N \l_abbr_count_tl

\NewDocumentCommand{\abbr}{m m O{#1} m m O{#4}}{
	\expandafter\newcommand\csname#3\endcsname[1][]{
		\seq_if_in:NnTF \g_abbrs {#1} {
			\prop_get:NnN \g_abbr_counts {#1} \l_abbr_count_tl
			\prop_gput:Nnx \g_abbr_counts {#1} {\int_eval:n {\l_abbr_count_tl + 1}}
			\hyperref[#1]{#1}
		} {
			\seq_gput_left:Nn \g_abbrs {#1}
			\prop_gput:Nnn \g_abbr_counts {#1} {1}
			\expandafter\gdef\csname#1@def\endcsname{#2}
			\phantomsection\label{#1}
			\str_if_eq:nnTF{##1}{}{\emph{#2}}{##1}~(\hyperref[#1]{#1})
		}
	}
	\expandafter\newcommand\csname#6\endcsname[1][]{
		\seq_if_in:NnTF \g_abbrs {#1} {
			\prop_get:NnN \g_abbr_counts {#1} \l_abbr_count_tl
			\prop_gput:Nnx \g_abbr_counts {#1} {\int_eval:n {\l_abbr_count_tl + 1}}
			\hyperref[#1]{#4}
		} {
			\expandafter\gdef\csname#1@def\endcsname{#5}
			\seq_gput_left:Nn \g_abbrs {#1}
			\prop_gput:Nnn \g_abbr_counts {#1} {1}
			\phantomsection\label{#1}
			\str_if_eq:nnTF{##1}{}{\emph{#5}}{##1}~(\hyperref[#1]{#4})
		}
	}
}

\ExplSyntaxOff
\makeatother

\abbr{ANN}{artificial neural network}{ANNs}{artificial neural networks}
\abbr{DL}{deep learning}{DLs}{deep learnings}
\abbr{SGD}{stochastic gradient descent}{SGDs}{stochastic gradient descent}
\abbr{Adam}{adaptive moment estimation \SGD}{Adams}{stochastic gradient descent}
\abbr{EMA}{exponential moving average}{EMAs}{exponential moving averages}
\abbr{SDE}{stochastic differential equation}{SDEs}{stochastic differential equations}
\abbr{ODE}{ordinary differential equation}{ODEs}{ordinary differential equations}
\abbr{GADAM}{geometrically weighted averaged \Adam}{GADAMs}{geometrically weighted averaged \Adam s}
\abbr{PADAM}{partially arithmetically averaged \Adam}{PADAMs}{partially arithmetically averaged \Adam s}
\abbr{PDE}{partial differential equation}{PDEs}{partial differential equations}
\abbr{GD}{gradient descent}{GDs}{gradient descent}
\abbr{DKM}{deep Kolmogorov method}{DKMs}{deep Kolmogorov methods}
\abbr{DK}{deep Kolmogorov}{DKs}{deep Kolmogorovs}
\abbr{AI}{artificial intelligence}{AIs}{artificial intelligences}
\abbr{PINN}{physics-informed neural network}{PINNs}{physics-informed neural networks}
\abbr{deep BSDE}{deep backward stochastic differential equation}[deepBSDE]{deep BSDEs}{deep backward stochastic differential equations}[deepBSDEs]
\abbr{DRM}{deep Ritz method}{DRMs}{deep Ritz methods}
\abbr{GELU}{Gaussian error linear unit}{GELUs}{Gaussian error linear units}
\abbr{ReLU}{rectified linear unit}{ReLUs}{rectified linear units}
\abbr{ResNet}{residual neural network}{ResNets}{residual neural networks}
\abbr{OC}{optimal control}{OCs}{optimal controls}
\abbr{HJB}{Hamiltonian--Jacobi--Bellman}{HJBs}{Hamiltonian--Jacobi--Bellmans}
\abbr{SOP}{stochastic optimization problem}{SOPs}{stochastic optimization problems}

\title{Averaged Adam accelerates stochastic optimization in\\ 
the training of deep neural network approximations for\\ 
partial differential equation and optimal control problems}

\author{Steffen Dereich$^{1}$, Arnulf Jentzen$^{2,3}$, and Adrian Riekert$^{4}$\bigskip\\
\small{$^1$ Institute for Mathematical Stochastics, University of M\"unster,}\vspace{-0.1cm}\\
\small{Germany; e-mail: \texttt{steffen.dereich}\textcircled{\texttt{a}}\texttt{uni-muenster.de}}\smallskip\\
\small{$^2$ School of Data Science and Shenzhen Research Institute of Big Data,}\vspace{-0.1cm}\\
\small{The Chinese University of Hong Kong, Shenzhen (CUHK-Shenzhen),}\vspace{-0.1cm}\\
\small{China; e-mail: \texttt{ajentzen}\textcircled{\texttt{a}}\texttt{cuhk.edu.cn}}\smallskip\\
\small{$^3$ Applied Mathematics: Institute for Analysis and Numerics,}\vspace{-0.1cm}\\
\small{University of M\"unster, Germany; e-mail: \texttt{ajentzen}\textcircled{\texttt{a}}\texttt{uni-muenster.de}}\smallskip\\
\small{$^4$ Applied Mathematics: Institute for Analysis and Numerics,}\vspace{-0.1cm}\\
\small{University of M\"unster, Germany; e-mail: \texttt{ariekert}\textcircled{\texttt{a}}\texttt{uni-muenster.de}}}

\date{\today}

\begin{document}
	
\maketitle
	
\begin{abstract}
Deep learning methods -- usually consisting of a class of 
deep neural networks (DNNs) trained by a stochastic gradient descent (SGD) 
optimization method -- are nowadays omnipresent in data-driven learning problems 
as well as in scientific computing tasks such as optimal control (OC) 
and partial differential equation (PDE) problems. 
In practically relevant learning tasks, often not the plain-vanilla 
standard SGD optimization method is employed to train the considered class 
of DNNs but instead more sophisticated adaptive and accelerated variants 
of the standard SGD method such as the popular Adam optimizer are used. 
Inspired by the classical Polyak--Ruppert averaging approach,  
in this work we apply averaged variants of the Adam optimizer 
to train DNNs to approximately solve exemplary scientific computing problems 
in the form of PDEs and OC problems. 
We test the averaged variants of Adam in a series of learning problems 
including physics-informed neural network (PINN), 
deep backward stochastic differential equation (deep BSDE), 
and deep Kolmogorov approximations for PDEs 
(such as heat, Black--Scholes, Burgers, and Allen--Cahn PDEs), 
including DNN approximations for OC problems, 
and including DNN approximations for image classification problems 
(ResNet for CIFAR-10). 
In each of the numerical examples the employed averaged variants of Adam 
outperform the standard Adam and the standard SGD optimizers, particularly, in the situation 
of the scientific machine learning problems. 
The {\sc Python} source codes for the numerical experiments associated to this work can be found on {\sc GitHub} at 
\href{https://github.com/deeplearningmethods/averaged-adam}{https://github.com/deeplearningmethods/averaged-adam}. 
\end{abstract}
	
\tableofcontents

\section{Introduction}

\DL[\emph{Deep learning}]\ methods are nowadays not only highly used 
to approximately solve data driven learning problems -- 
such as those occurring in \AI\ based chatbot systems 
(cf., \eg, \cite{Brownlanguagemodel2020,Yihengsummary2023}) and \AI\ based text-to-image creation models 
(cf., \eg, \cite{Rombachetal2022arXiv,Sahariaimagen2022_IMAGEN,ramesh2021zeroshottexttoimagegeneration_DALLE}) -- but \DL\ methods are these days also extensively used to 
approximate solutions of scientific computing problems 
such as \PDEs\ and \OC\ problems 
(cf., \eg, the overview articles \cite{beck2020overview,Weinan2021,Germainetal2021arXivOverviewArticle,MR4795589,MR4457972}).

\DL\ methods typically consist of a class of 
deep \ANNs\ that are trained by a \SGD\ optimization method. 
In practically relevant learning problems, often not the plain-vanilla standard \SGD\ optimization method 
is employed to train the considered class of deep \ANNs\ but instead more sophisticated 
adaptive and accelerated variants of the standard \SGD\ method are used 
(cf., \eg, \cite{bach2024learning,Ruder2017overview,JentzenBookDeepLearning2023,Sun2019} 
for overviews and monographs). 
Maybe the most popular variant of such adaptive and accelerated \SGD\ methods 
is the famous \Adam\ optimizer proposed in 2014 by Kingma \& Ba 
(see \cite{KingmaBa2014}).

Inspired by the classical Polyak--Ruppert averaging approach \cite{Ruppert1988,MR1071220} (cf.\ also \cite{MR1167814}), 
several averaged variants of \SGD\ optimization methods have been considered 
in the literature as well. 
In particular, we refer, \eg, to \cite{morales-brotons2024exponential,SandlerZhmoginov2023,Defazioetal2024arXiv,ZhangChoromanskaLeCun2014arXiv,BusbridgeRamapuram2023EMA,GuoJin2023SWA,IzmailovPodoprikhin2018SWA,Athiwaratkun2018Average} 
and the references therein for works proposing and testing 
\SGD\ methods involving suitable averaging techniques 
and we refer, \eg, to \cite{AhnCutkosky2024arXiv,MR4184368,DereichMuller_Gronbach2019,MR4580893,GadatPanloup2017arXiv,AhnMagakyanCutkosky2024arXiv,MandtHoffmanBlei2017arXiv} 
and the references therein for articles studying averaged variants of \SGD\ methods analytically.

In this work we apply different averaged variants of the \Adam\ optimizer (see \cref{sec:numerics}) 
to train deep \ANNs\ to approximately solve exemplary scientific computing 
and image classification problems.  
%
%
%
Specifically, we study the considered averaged variants of the \Adam\ optimizer numerically in a series of learning problems 
\begin{itemize}
\item 
including polynomial regression problems
(see \cref{subsec:poly_reg}), 

\item 
including deep \ANN\ approximations 
for  explicitly given high-dimensional target functions
(see \cref{subsec:simple_supervised}), 

\item 
including 
\PINN\ (see \cref{ssub:PINNs}), 
\deepBSDE\ (see \cref{ssub:deepBSDE}), 
and \DK\ (see \cref{ssub:DKM}) approximations 
for \PDEs\ (such as heat, Black--Scholes, Burgers, Allen--Cahn, and \HJB\ \PDEs), 

\item 
including deep \ANN\ approximations for stochastic \OC\ problems 
(see \cref{ssub:OC}), 
and

\item 
including residual deep \ANN\ approximations 
for the CIFAR-10 image classification dataset 
(see \cref{ssub:CIFAR-10}). 
\end{itemize}
In each of the considered numerical examples 
the suggested optimizers 
outperform the standard \Adam\ and the standard \SGD\ optimizers, 
particularly, in the situation of 
the considered scientific computing problems.
Taking this into account, 
we strongly suggest to further study and employ averaged variants of the Adam optimizers 
when solving \PDE, \OC, or related scientific computing problems by 
means of deep learning approximation methods. 
The {\sc Python} source codes for each of the performed numerical simulations can 
be found on {\sc GitHub} at \href{https://github.com/deeplearningmethods/averaged-adam}{https://github.com/deeplearningmethods/averaged-adam}.

%
%
%
%

\subsubsection*{Structure of this article}

The remainder of this work is organized as follows. 
In \cref{sec:methods} we recall the concept of the standard \Adam\ optimizer 
and we describe in detail the specific averaged variants of the \Adam\ optimizer 
that we employ in our numerical simulations. 
In \cref{sec:numerics} we apply the considered averaged variants of the \Adam\ optimizer 
to several scientific computing and image classification problems 
and compare the obtained approximation errors with those 
of the standard \SGD\ and the standard \Adam\ optimizers. 
Finally, in \cref{sec:conclusion} we briefly summarize the findings of this work 
and also outline directions of future research.

\section{Averaged Adam optimizers}
\label{sec:methods}

\subsection{Standard Adam optimizer}

For convenience of the reader we recall within this subsection 
in \cref{def:adam} and \cref{alg:adam} below the concept of 
the ``standard'' \Adam\ optimizer from Kingma \& Ba~\cite{KingmaBa2014}. 
\cref{def:adam} is a modified variant of \cite[Definition~7.9.1]{JentzenBookDeepLearning2023}.


\begin{definition}[Standard \Adam\ optimizer]
\label{def:adam}
Let 
$ \defaultParamDim, \defaultdataDim \in \N $, 
$
  ( \gamma_n )_{ n \in \N } 
  \allowbreak 
  \subseteq \R
$, 
$
  ( J_n )_{ n \in \N } \subseteq \N 
$,
$ 
  ( \alpha_n )_{ n \in \N } \subseteq [0,1)
$, 
$
  ( \beta_n )_{ n \in \N } \subseteq [0,1)
$, 
$
  \varepsilon \in (0,\infty) 
$, 
let $ ( \Omega, \mathcal{F}, \P ) $ be a probability space, 
for every $ n, j \in \N $ let
$
  X_{ n, j } \colon \Omega \to \R^{ \defaultdataDim }
$
be a random variable, 
let 
$
  \defaultStochLoss \colon \R^\defaultParamDim \times \R^{ \defaultdataDim } \to \R
$ 
be differentiable, 
let 
$
  \defaultStochGradient = (\defaultStochGradient_1,\ldots, \defaultStochGradient_\defaultParamDim)  \colon \R^\defaultParamDim \times \R^{ \defaultdataDim } 
   \to \R^{ \defaultParamDim }
$ 
satisfy for all 
$ \theta \in \R^{ \defaultParamDim } $,
$ x \in \R^{ \defaultdataDim } $
that
\begin{equation}
\label{eq:Adam_generalized_gradient}
  \defaultStochGradient(\theta,x) = \nabla_{ \theta } \defaultStochLoss( \theta, x )
  ,
\end{equation}
and let 
$ 
  \Theta = (\Theta^{(1)}, \ldots, \Theta^{(\defaultParamDim)}) \colon \N_0 \times \Omega \to \R^\defaultParamDim
$
be a function. 
Then we say that $ \Theta $ is the 
\Adam\ process 
for $ \defaultStochLoss $
with hyperparameters 
$ ( \alpha_n )_{ n \in \N } $, 
$ ( \beta_n )_{ n \in \N } $, 
$ ( \gamma_n )_{ n \in \N } $, 
$ \varepsilon \in (0,\infty) $, 
batch-sizes 
$ ( J_n )_{ n \in \N } $, 
initial value $ \Theta_0 $, and data $ ( X_{ n, j } )_{ (n,j) \in \N^2 } $ 
if and only if there exist  
$ 
  \momentprocess = (\momentprocess^{(1)},\ldots,\momentprocess^{(\defaultParamDim)}) \colon \N_0\times \Omega \to \R^\defaultParamDim
$
and 
$ 
  {\bf v} = ({\bf v}^{(1)},\ldots,{\bf v}^{(\defaultParamDim)}) \colon \N_0\times \Omega \to \R^\defaultParamDim
$ 
such that for all $ n \in \N $, $ i \in \{ 1, 2, \dots, \defaultParamDim \} $ 
it holds that 
\begin{equation}
  \momentprocess_0 = 0, \qquad 
  \momentprocess_n = \alpha_n \, \momentprocess_{n-1} + (1-\alpha_n)\br*{\frac{1}{J_n}\sum_{j = 1}^{J_n}
  \defaultStochGradient( \Theta_{n-1}, X_{n,j} )},
\end{equation}
\begin{equation}
  {\bf v}_0 = 0, 
  \qquad 
  {\bf v}_n^{(i)} = \beta_n\,{\bf v}_{n-1}^{(i)} + (1-\beta_n)\br*{\frac{1}{J_n} 
  \sum_{j = 1}^{J_n} \defaultStochGradient_{i}(\Theta_{n-1}, X_{n,j})}^2, 
\end{equation}
\begin{equation}
  \andq 
  \Theta_n^{(i)} = \Theta_{n-1}^{(i)} 
  - 
  \gamma_n 
  \,
  {\textstyle 
    \left[
      \varepsilon 
      + 
      \left[ 
        \frac{ 
          {\bf v}_n^{ (i) } 
        }{ 
          ( 1 - \prod_{ k = 1 }^n \beta_k ) 
        } 
      \right]^{ \nicefrac{1}{2} } 
    \right]^{ - 1 } 
  } 
  \br*{
    \frac{ \momentprocess_n^{ (i) } }{ ( 1 - \prod_{ k = 1 }^n \alpha_k ) }
  } .
\end{equation}
\end{definition}

In \cref{def:adam} 
the hyperparameter $ \gamma $ can be referred to as \emph{learning rate}, 
the hyperparameter $ \alpha $ can be referred to as \emph{momentum decay factor}, 
the hyperparameter $ \beta $ can be referred to as \emph{second moment decay factor}, 
and the hyperparameter $ \varepsilon $ can be referred to as 
\emph{regularizing parameter} 
(cf., \eg, \cite[Definition~7.9.1]{JentzenBookDeepLearning2023}). 
In {\sc PyTorch} the hyperparameters $ ( \alpha_n )_{ n \in \N } $, 
$ ( \beta_n )_{ n \in \N } $, $ ( \gamma_n )_{ n \in \N } $, and $ \varepsilon $ 
for the \Adam\ optimizer in \cref{def:adam} 
are by default chosen to satisfy 
for all $ n \in \N $ that 
$
  \alpha_n = 0.9
$, 
$ 
  \beta_n = 0.999 
$, 
$
  \gamma_n = 0.001
$, 
and 
$
  \varepsilon = 10^{ - 8 }
$
(cf., \eg, \cite{PyTorchAdamWeb2025,KingmaBa2014}).
In \cref{alg:adam} below we present a pseudo-code 
for the standard \Adam\ optimizer in \cref{def:adam}.

\begin{myalgorithm}{Standard \Adam\ optimizer}{alg:adam}
\begin{algorithmic}[1]
\Statex\textbf{Setting:} The mathematical objects introduced in \cref{def:adam}
\Statex\textbf{Input:} 
$ \numberofsteps \in \N $
\Statex\textbf{Output:} 
\Adam\ process $ \Theta_{ \numberofsteps } \in \R^\defaultParamDim $ after $ \numberofsteps $ steps 

\vspace{-2mm} 
\noindent\hspace*{-8.1mm}\rule{\dimexpr\linewidth+9.24mm}{0.4pt}  
\vspace{-3mm} 
\State $ \theta \gets \Theta_0 $
\State $ \momentprocess \gets 0 $
\State $ \Momentprocess \gets 0 $
\For{$ n \in \cu{ 1, 2, \dots, \numberofsteps } $}
  \State $ g \gets ( J_n )^{ - 1 } \sum_{ j = 1 }^{ J_n } \defaultStochGradient ( \theta , X_{ n, j } ) $
  \State $ \momentprocess \gets \alpha_n \momentprocess + ( 1 - \alpha_n ) g $
  \State $ \Momentprocess \gets \beta_n \Momentprocess + ( 1 - \beta_n ) g^{ \otimes 2 }$ \mycomment{Square $ g^{ \otimes 2 } $ is understood componentwise}
  \State $ \hat{\momentprocess} \gets \momentprocess / (1 - \prod_{ k = 1 }^n \alpha_k ) $
  \State $ \hat{\Momentprocess} \gets \Momentprocess / (1 - \prod_{ k = 1 }^n \beta_k ) $
  \State $ \theta \gets \theta - \gamma_n \hat{\momentprocess} \otimes ( \hat{\Momentprocess}^{ \otimes (1 / 2) } + \varepsilon)^{ \otimes ( - 1 ) } $ 
  \mycomment{Root $ \hat{\Momentprocess}^{ \otimes (1 / 2) } $ is understood componentwise}
\EndFor
\noindent 
\State \Return $\theta$
\end{algorithmic}
\end{myalgorithm}

We refer, \eg, to \cite{Barakat_2021_cvg,li2023convergenceadamrelaxedassumptions,DereichJentzen2024arXiv_Adam,ReddiKale2019,Defossez2022} and the references therein 
for error and convergence rate analyses for the \Adam\ optimizer.

\subsection{General averaged Adam optimizers}
\label{ssec:averaged_adam}

In order to be in the position to precisely describe the optimization methods that we employ in the numerical simulations 
in \cref{sec:numerics} below, we present in this subsection (see \cref{def:averaged_adam_general} below) 
for the convenience of the reader a general class of averaged variants of the \Adam\ optimizer.

\cfclear
\begin{definition}[General averaged \Adam\ optimizer]
\label{def:averaged_adam_general}
Let 
$ \defaultParamDim, \defaultdataDim \in \N $, 
$
  ( \gamma_n )_{ n \in \N } 
  \allowbreak 
  \subseteq \R
$, 
$
  ( 
    J_n 
  )_{ n \in \N } 
  \allowbreak 
  \subseteq \N 
$,
$ 
  ( \alpha_n )_{ n \in \N } \subseteq [0,1)
$, 
$
  ( \beta_n )_{ n \in \N } \subseteq [0,1)
$, 
$
  ( \delta_{ n, m } )_{ (n, m) \in ( \N_0 )^2 } \subseteq \R
$, 
$
  \varepsilon \in (0,\infty) 
$, 
let $ ( \Omega, \mathcal{F}, \P ) $ be a probability space, 
for every $ n, j \in \N $ let
$
  X_{ n, j } \colon \Omega \to \R^{ \defaultdataDim }
$
be a random variable, 
let 
$
  \defaultStochLoss \colon \R^\defaultParamDim \times \R^{ \defaultdataDim } \to \R
$ 
be differentiable, 
and let 
$ 
  \Theta 
  \colon \N_0 \times \Omega \to \R^\defaultParamDim
$
be a function. 
Then we say that 
$ \Theta $ is the averaged \Adam\ process 
with loss $ \defaultStochLoss $, 
learning rates $(\gamma_n)_{n \in \N}$, batch sizes $(J_n)_{n\in \N}$, 
momentum decay factors $(\alpha_n)_{n \in \N}$, 
second moment decay factors $(\beta_n)_{n \in \N}$, 
regularizing factor $ \varepsilon \in (0,\infty) $, 
initial value $ \Theta_0 $, 
data $ ( X_{n,j} )_{ (n,j) \in \N^2 } $, 
and 
averaging weights $ ( \delta_{ n, m } )_{ (n,m) \in ( \N_0 )^2 } $
if and only if there exists   
$ 
  \vartheta \colon \N_0 \times \Omega \to \R^\defaultParamDim
$
such that 
\begin{enumerate}[label=(\roman*)]
\item 
it holds that $ \vartheta $ is the \Adam\ process 
with loss $ \defaultStochLoss $, 
learning rates $ ( \gamma_n )_{ n \in \N } $, 
batch sizes $ ( J_n )_{ n \in \N } $, 
momentum decay factors $ ( \alpha_n )_{ n \in \N } $, 
second moment decay factors $ ( \beta_n )_{ n \in \N } $, 
regularizing factor $ \varepsilon \in (0,\infty) $, 
initial value $ \Theta_0 $, and data $ ( X_{ n, j } )_{ (n,j) \in \N^2 } $ 
\cfadd{def:adam}\cfload\ and 
\item 
it holds for all $ n \in \N $ that 
\begin{equation}
\label{eq:general_averaged_adam}
  \Theta_n = 
  \sum_{ m = 0 }^n \delta_{ n, m } \vartheta_m
  .
\end{equation}
\end{enumerate}
\end{definition}

In the following list we present a few 
special cases of \cref{def:averaged_adam_general} 
by choosing specific values for the family 
of averaging weights 
$
  \delta_{ n, m } \in \R
$, 
$
  (n,m) \in ( \N_0 )^2
$,
in \cref{eq:general_averaged_adam}. 
\begin{enumerate}[label=(\Roman*)]

\item 
\emph{Standard \Adam}: 
Consider \cref{def:averaged_adam_general} 
and assume for all $ n, m \in \N_0 $ 
that $ \delta_{ n, m } = \mathbbm{1}_{ \{ n \} }( m ) $. 
Then it holds for all $ n \in \N_0 $ that 
\begin{equation}
  \Theta_n = \vartheta_n 
\end{equation}
and, in this situation, 
the averaged \Adam\ process 
in \cref{def:averaged_adam_general} 
reduces to the standard \Adam\ process in \cref{def:adam}.

\item 
\label{item:partially_arithmetically}
\emph{Arithmetic average of \Adam\ over the last $ \AvgParam $ steps}: 
Consider \cref{def:averaged_adam_general}, 
let $ \AvgParam \in \N $, and 
assume for all 
$ n, m \in \N_0 $ that
$ 
  \delta_{ n, m } 
  = ( \AvgParam + 1 )^{ - 1 } \mathbbm{1}_{ [ 0, \AvgParam ] }( n - m )
$. 
Then it holds for all $ n \in \N \cap [ \AvgParam, \infty ) $ that 
\begin{equation}
\label{eq:adam_avg_last_A_steps}
  \Theta_n 
  = \frac{ \sum_{ k = n + 1 - \AvgParam }^n \vartheta_k }{ \AvgParam }
  .
\end{equation}
The choice in \cref{eq:adam_avg_last_A_steps} 
is the subject of \cref{def:averaged_adam_arithmetic}, \cref{alg:adam_avg_a}, 
and \cref{alg:adam_avg_a2} in \cref{subsec:arithmetic}. 
Moreover, in the case where $ \AvgParam = 999 $ 
in \cref{eq:adam_avg_last_A_steps} 
we present in \cref{sec:numerics}
a series of numerical simulations 
for this averaged variant of the \Adam\ optimizer. 
We also refer to this type of averaging of the \Adam\ optimizer 
as \emph{partially arithmetically averaged \Adam\ optimizer}.

\item 
\label{item:geometrically}
\emph{Geometrically weighted averages of \Adam}: 
Consider \cref{def:averaged_adam_general}, 
let $ ( \varrho_n )_{ n \in \N } \subseteq \R $ 
satisfy 
for all $ n \in \N $ 
that 
$
  \varrho_n = 1 - \delta_{ n, n } 
$, 
and assume for all $ n, m \in \N_0 $ 
with $ m < n $ 
that 
$
  \delta_{ n, m } = 
  ( 1 - \delta_{ n, n } ) 
  \delta_{ n - 1, m }
$.
Then 
it holds for all $ n \in \N $ that 
\begin{equation}
\label{eq:adam_avg_geo}
\begin{split}
  \Theta_n 
& 
  = 
  \sum_{ k = 0 }^n 
  \delta_{ n, k }
  \vartheta_k
  =
  \left[ 
    \sum_{ k = 0 }^{ n - 1 } 
    \delta_{ n, k }
    \vartheta_k
  \right] 
  +
  \delta_{ n, n }
  \vartheta_n
  =
  \left[ 
    \sum_{ k = 0 }^{ n - 1 } 
    ( 1 - \delta_{ n, n } )
    \delta_{ n - 1, k }
    \vartheta_k
  \right] 
  +
  \delta_{ n, n }
  \vartheta_n
\\ &  
  =
  ( 1 - \delta_{ n, n } )
  \Theta_{ n - 1 }
  +
  \delta_{ n, n }
  \vartheta_n 
  =
  \varrho_n \Theta_{ n - 1 } + ( 1 - \varrho_n ) \vartheta_n 
  .
\end{split}
\end{equation}
\Hence for all $ n \in \N $ that 
\begin{equation}
  \Theta_n 
  = 
  \left[ 
    \prod_{ k = 1 }^n
    \varrho_k 
  \right] 
  \Theta_0
  +
  \sum_{ k = 1 }^n
  \left(
  \left[ 
    \prod_{ v = k + 1 }^{ n } \varrho_v
  \right] 
  \left( 1 - \varrho_k \right)
  \vartheta_k
  \right)
\end{equation}
The choice in \cref{eq:adam_avg_geo} 
is the subject of \cref{def:averaged_adam_geometric} 
and \cref{alg:adam_avg_g} in \cref{ssub:geometrically}. 
Moreover, in the case where it holds for all $ n \in \N $ that 
$ \varrho_n = 1 - \delta_{ n, n } = 0.999 $ 
in \cref{eq:adam_avg_geo} we present in \cref{sec:numerics}
a series of numerical simulations 
for this averaged variant of the \Adam\ optimizer. 
In the scientific literature the type of averaging in \cref{eq:adam_avg_geo} 
is referred to as \EMA\ (cf., \eg, 
\cite{morales-brotons2024exponential,SandlerZhmoginov2023,AhnCutkosky2024arXiv,BusbridgeRamapuram2023EMA,IzmailovPodoprikhin2018SWA,Athiwaratkun2018Average,GuoJin2023SWA}).

%

\item 
\emph{Arithmetic average of \Adam\ over all steps since the $ \AvgParam $-th step}: 
Consider \cref{def:averaged_adam_general}, 
let $ \AvgParam \in \N_0 $, 
and assume for all $ n, m \in \N_0 $ that
$
  \delta_{ n, m } = 
  ( n + 1 - \AvgParam )^{ - 1 }
  \mathbbm{1}_{ [ \AvgParam, n ] }( m )
$.
Then it holds for all $ n \in \N \cap [ \AvgParam, \infty ) $ that 
\begin{equation}
  \Theta_n 
  = \frac{ \sum_{ k = \AvgParam }^n \vartheta_k }{ n + 1 - \AvgParam }
\end{equation}
\Hence for all $ n \in \N \cap ( \AvgParam, \infty ) $ that 
\begin{equation}
\label{eq:recursion_avg_A}
\begin{split}
  \Theta_n
&
  = 
  \frac{ \sum_{ k = \AvgParam }^n \varTheta_k }{ n + 1 - \AvgParam }
  =
  \frac{ \sum_{ k = \AvgParam }^{ n - 1 } \varTheta_k }{ n + 1 - \AvgParam }
  +
  \frac{
    \vartheta_n 
  }{
    n + 1 - \AvgParam
  }
\\ & 
  =
  \Theta_{ n - 1 } 
  \left[ 
    \frac{ ( n - 1 ) + 1 - \AvgParam }{ n + 1 - \AvgParam }
  \right] 
  +
  \frac{
    \vartheta_n 
  }{
    n + 1 - \AvgParam
  }
\\ &
  =
  \left[ 
    \frac{ n - \AvgParam }{ n + 1 - \AvgParam }
  \right] 
  \Theta_{ n - 1 } 
  +
  \left(
    1 
    -
    \left[
      \frac{ n - \AvgParam }{ n + 1 - \AvgParam }
	\right]
  \right) 
  \vartheta_n 
  .
\end{split}
\end{equation}

\item 
\emph{Arithmetic average of \Adam\ over all previous steps}: 
Consider \cref{def:averaged_adam_general} 
and assume for all $ n, m \in \N_0 $ that
$
  \delta_{ n, m } = 
  ( n + 1 )^{ - 1 }
$.
Then it holds for all $ n \in \N $ that 
\begin{equation}
  \Theta_n 
  = \frac{ \sum_{ k = 0 }^n \vartheta_k }{ n + 1 }
\end{equation}
Combining this and \cref{eq:recursion_avg_A} 
proves for all $ n \in \N $ that 
\begin{equation}
\begin{split}
  \Theta_n
  =
  \left[ 
    \frac{ n }{ n + 1 }
  \right] 
  \Theta_{ n - 1 } 
  +
  \left(
    1 
    -
    \left[
      \frac{ n }{ n + 1 }
    \right]
  \right) 
  \vartheta_n 
  .
\end{split}
\end{equation}
This type of averaging corresponds to the classical Polyak--Ruppert 
averaging approach (see \cite{Ruppert1988,MR1071220,MR1167814}).

\end{enumerate}

\subsection{Partially arithmetically averaged Adam optimizer}
\label{subsec:arithmetic}


As we employ the partially arithmetically averaged variant of the \Adam\ optimizer 
in \cref{item:partially_arithmetically} in \cref{ssec:averaged_adam} above 
in each of our numerical simulations in \cref{sec:numerics}, 
we describe this type of averaging of \Adam\ and its specific implementations 
in \cref{def:averaged_adam_arithmetic}, \cref{alg:adam_avg_a}, and \cref{alg:adam_avg_a2} 
within this subsection in more details.

\begin{definition}[Partially arithmetically averaged \Adam\ optimizer]
\label{def:averaged_adam_arithmetic}
Let 
$ 
  \defaultParamDim, \defaultdataDim, \AvgParam \in \N 
$, 
$
  \allowbreak 
  ( 
    \gamma_n 
    \allowbreak 
  )_{ n \in \N } 
  \allowbreak 
  \subseteq \R
$, 
$
  ( 
    J_n 
  )_{ n \in \N } 
  \allowbreak 
  \subseteq \N 
$,
$ 
  ( \alpha_n )_{ n \in \N } \subseteq [0,1)
$, 
$
  ( \beta_n )_{ n \in \N } \subseteq [0,1)
$,
$
  \varepsilon \in (0,\infty) 
$, 
let $ ( \Omega, \mathcal{F}, \P ) $ be a probability space, 
for every $ n, j \in \N $ let
$
  X_{ n, j } \colon \Omega \to \R^{ \defaultdataDim }
$
be a random variable, 
let 
$
  \defaultStochLoss \colon \R^\defaultParamDim \times \R^{ \defaultdataDim } \to \R
$ 
be differentiable, 
and let 
$ 
  \Theta 
  \colon \N_0 \times \Omega \to \R^\defaultParamDim
$
be a function. 
Then we say that 
$ \Theta $ is the $ \AvgParam $-partially averaged \Adam\ process 
with loss $ \defaultStochLoss $, 
learning rates $(\gamma_n)_{n \in \N}$, batch sizes $(J_n)_{n\in \N}$, 
momentum decay factors $(\alpha_n)_{n \in \N}$, 
second moment decay factors $(\beta_n)_{n \in \N}$, 
regularizing factor $ \varepsilon \in (0,\infty) $, 
initial value $ \Theta_0 $, 
and 
data $ (X_{n,j})_{ (n,j) \in \N^2 } $
if and only if there exists   
$ 
  \vartheta \colon \N_0 \times \Omega \to \R^\defaultParamDim
$
such that 
\begin{enumerate}[label=(\roman*)]
\item 
it holds that $ \vartheta $ is the \Adam\ process 
with loss $ \defaultStochLoss $, 
learning rates $ ( \gamma_n )_{ n \in \N } $, 
batch sizes $ ( J_n )_{ n \in \N } $, 
momentum decay factors $ ( \alpha_n )_{ n \in \N } $, 
second moment decay factors $ ( \beta_n )_{ n \in \N } $, 
regularizing factor $ \varepsilon \in (0,\infty) $, 
initial value $ \Theta_0 $, and data $ ( X_{ n, j } )_{ (n,j) \in \N^2 } $ 
and 
\item 
it holds for all $ n \in \N \cap [ \AvgParam, \infty ) $ that 
\begin{equation}
  \Theta_n 
  = 
  \frac{ 1 }{ \AvgParam }
\textstyle 
  \left[ 
    \sum\limits_{ k = n + 1 - \AvgParam }^n \vartheta_k
  \right] 
  .
\end{equation}
\end{enumerate}
\end{definition}

In the following we describe 
in \cref{alg:adam_avg_a,alg:adam_avg_a2}
two different concrete implementations 
of the method in \cref{def:averaged_adam_arithmetic} above.

\begin{myalgorithm}{\Adam\ with partial arithmetic averaging}{alg:adam_avg_a}
\begin{algorithmic}[1]
\Statex\textbf{Setting:} The mathematical objects introduced in \cref{def:averaged_adam_arithmetic}
\Statex\textbf{Input:} 
$ \numberofsteps \in \N $
\Statex\textbf{Output:} 
$ \AvgParam $-partially averaged \Adam\ process $ \Theta_{ \numberofsteps } \in \R^\defaultParamDim $ after $ \numberofsteps $ steps 

\vspace{-2mm} 
\noindent\hspace*{-8.1mm}\rule{\dimexpr\linewidth+9.24mm}{0.4pt}  
\vspace{-3mm} 

\State $ \vartheta \gets \Theta_0 $
\State $ \theta \gets \Theta_0 $
\State $ \phi_0 \gets \Theta_0 $
\State $ \momentprocess \gets 0 $
\State $ \Momentprocess \gets 0 $
\For{$ \nos \in \cu{ 1, 2, \dots, \numberofsteps } $}
\State $ g \gets ( J_n )^{ - 1 } \sum_{ j = 1 }^{ J_n } \defaultStochGradient( \vartheta , X_{ n, j } ) $
\State $ \momentprocess \gets \alpha_n \momentprocess + ( 1 - \alpha_n ) g $
\State $ \Momentprocess \gets \beta_n \Momentprocess + ( 1 - \beta_n ) g^{ \otimes 2 } $ \mycomment{Square $ g^{ \otimes 2 } $ is understood componentwise}
\State $ \hat{\momentprocess} \gets \momentprocess / (1 - \prod_{k=1}^n \alpha_k )$
\State $ \hat{\Momentprocess} \gets \Momentprocess / (1 - \prod_{k=1}^n \beta_k )$
\State $ \vartheta \gets \vartheta - \gamma_n \hat{\momentprocess} \otimes ( \hat{\Momentprocess}^{ \otimes (1 / 2) } + \varepsilon)^{ \otimes ( - 1 ) } 
$ \mycomment{Root $ \hat{\Momentprocess}^{ \otimes (1/2) } $ is understood componentwise}
\State $ \theta \gets \theta + \AvgParam^{-1} ( \vartheta - \phi_0 ) $ \mycomment{Update averaged iterate}
\If{$ n < \AvgParam $ }
\State $\phi_i \gets \vartheta$
\Else 
\State $(\phi_0, \phi_1, \ldots, \phi_{N - 1 } ) \gets (\phi_1, \phi_2, \ldots, \phi_{N-1}, \vartheta)$ \mycomment{Store previous $N$ iterates}
\EndIf
\EndFor
\noindent 
\State \Return $ \theta $
\end{algorithmic}
\end{myalgorithm}

\cref{alg:adam_avg_a} has the disadvantage that one needs to store all previous $ \AvgParam $ iterates.
In the following pseudocode in \cref{alg:adam_avg_a2} below 
we decompose $ \AvgParam = \AvgPb \AvgPc $ and only update the averages every $ \AvgPc $ steps.
In this way one only needs to store the average over groups of $ \AvgPc $ iterates, 
i.e., $ \AvgPb $ instead of $ \AvgPb \AvgPc $ additional parameter vectors.
In our numerical simulations 
in \cref{sec:numerics} below we implement the method in \cref{def:averaged_adam_arithmetic} 
using \cref{alg:adam_avg_a2}, \eg, with the choice $ \AvgPb = 1 $, $ \AvgPc = 1000 $, $ \AvgParam = 1000 $ 
in \cref{alg:adam_avg_a2}.

\begin{myalgorithm}{\Adam\ with partial arithmetic averaging over groups}{alg:adam_avg_a2}
\begin{algorithmic}[1]
		\Statex\textbf{Setting:} The mathematical objects introduced in \cref{def:averaged_adam_arithmetic}
		\Statex\textbf{Input:} $ \AvgPb, \AvgPc, \numberofsteps \in \N $
		
		\Statex\textbf{Output:} $ \AvgParam $-partially averaged \Adam\ process $ \Theta_{ \numberofsteps } \in \R^\defaultParamDim $ after $ \numberofsteps $ steps 
		
		\vspace{-2mm} 
		\noindent\hspace*{-8.1mm}\rule{\dimexpr\linewidth+9.24mm}{0.4pt}  
		\vspace{-3mm} 
		\State $ \vartheta \gets \Theta_0 $
		\State $ \theta \gets \Theta_0 $
		\State $ \phi_0 \gets \Theta_0 $
		\State $ \AvgPx \gets 0 $
		\State $ \momentprocess \gets 0$
		\State $ \Momentprocess \gets 0 $
		\For{$ \nos \in \cu{ 1, 2, \dots, \numberofsteps } $}
		\State $ g \gets ( J_{ \nos } )^{ - 1 } \sum_{ j = 1 }^{ J_{ \nos } } \defaultStochGradient( \theta, X_{ \nos, j } ) $
		\State $ \momentprocess \gets \alpha_{ \nos } \momentprocess + ( 1 - \alpha_n ) g $
		\State $ \Momentprocess \gets \beta_{ \nos } \Momentprocess + ( 1 - \beta_n ) g^{ \otimes 2 } $ 
		\mycomment{Square $ g^{ \otimes 2 } $ is understood componentwise}
		\State $ \hat{\momentprocess} \gets \momentprocess / (1 - \prod_{ k = 1 }^n \alpha_k ) $
		\State $ \hat{\Momentprocess} \gets \Momentprocess / (1 - \prod_{ k = 1 }^n \beta_k ) $
		\State $ \vartheta \gets \vartheta - \gamma_{ \nos } \hat{\momentprocess} \otimes ( \hat{\Momentprocess}^{ \otimes (1/2) } + \varepsilon )^{ \otimes (-1) } $ 
		\mycomment{Root $ \Momentprocess^{ \otimes ( 1/2 ) } $ is understood componentwise}
		\State $ \AvgPx \gets \AvgPx + \AvgPc^{ - 1 } \vartheta$ 
		\mycomment{Update average of group of $ \AvgPc $ iterates}
		\If{$ \nos \equiv 0 \pmod{\AvgPc} $}
		\State $ \theta \gets \theta + \AvgPb^{ - 1 } ( \AvgPx - \phi_0 ) $
		\mycomment{Update averaged iterate}
		\If{$ \nos < \AvgPb \AvgPc $}
		\State $ \phi_{ \nos } \gets \AvgPx $
		\Else
		\mycomment{Store previous $ \AvgPb $ group averages} 
		\State $ (\phi_0, \phi_1, \dots, \phi_{ \AvgPb - 1 } ) \gets ( \phi_1, \phi_2, \dots, \phi_{ \AvgPb - 1 }, \AvgPx ) $ 
		\EndIf
		\State $ \AvgPx \gets 0 $
		\EndIf
		\EndFor
		
		\noindent 
		\State \Return $ \theta $
	\end{algorithmic}
\end{myalgorithm}

\subsection{Geometrically weighted averaged Adam optimizer}
\label{ssub:geometrically}

As we also use the 
geometrically weighted averaged variant of \Adam\ in \cref{item:geometrically} in \cref{ssec:averaged_adam} above 
in each of our numerical simulations in \cref{sec:numerics} below, 
we describe this type of averaging of \Adam\ and its implementation
in \cref{def:averaged_adam_geometric} and \cref{alg:adam_avg_g} 
within this subsection in more details.

\begin{definition}[Geometrically weighted averaged \Adam\ optimizer]
\label{def:averaged_adam_geometric}
Let 
$ \defaultParamDim, \defaultdataDim \in \N $, 
$
  ( \gamma_n )_{ n \in \N } 
  \allowbreak 
  \subseteq \R
$, 
$
  ( 
    J_n 
  )_{ n \in \N } 
  \allowbreak 
  \subseteq \N 
$,
$ 
  ( \alpha_n )_{ n \in \N } \subseteq [0,1)
$, 
$
  ( \beta_n )_{ n \in \N } \subseteq [0,1)
$, 
$
  ( \avgPa_n )_{ n \in \N } \subseteq \R
$, 
$
  \varepsilon \in (0,\infty) 
$, 
let $ ( \Omega, \mathcal{F}, \P ) $ be a probability space, 
for every $ n, j \in \N $ let
$
  X_{ n, j } \colon \Omega \to \R^{ \defaultdataDim }
$
be a random variable, 
let 
$
  \defaultStochLoss \colon \R^\defaultParamDim \times \R^{ \defaultdataDim } \to \R
$ 
be differentiable, 
let 
$
  \defaultStochGradient = (\defaultStochGradient_1,\ldots, \defaultStochGradient_\defaultParamDim)  \colon \R^\defaultParamDim \times \R^{ \defaultdataDim } 
   \to \R^{ \defaultParamDim }
$ 
satisfy for all 
$ \theta \in \R^{ \defaultParamDim } $,
$ x \in \R^{ \defaultdataDim } $
that
\begin{equation}
\label{eq:Adam_generalized_gradient3}
  \defaultStochGradient(\theta,x) = \nabla_{ \theta } \defaultStochLoss( \theta, x )
  ,
\end{equation}
and let 
$ 
  \Theta 
  \colon \N_0 \times \Omega \to \R^\defaultParamDim
$
be a function. 
Then we say that 
$ \Theta $ is the geometrically weighted averaged \Adam\ process 
with loss $ \defaultStochLoss $, 
learning rates $(\gamma_n)_{n \in \N}$, batch sizes $(J_n)_{n\in \N}$, 
momentum decay factors $(\alpha_n)_{n \in \N}$, 
second moment decay factors $(\beta_n)_{n \in \N}$, 
regularizing factor $ \varepsilon \in (0,\infty) $, 
initial value $ \Theta_0 $, 
data $ (X_{n,j})_{ (n,j) \in \N^2 } $, 
and 
averaging weights $ ( \avgPa_{ \nos } )_{ \nos \in \N } $
if and only if there exists   
$ 
  \vartheta \colon \N_0 \times \Omega \to \R^\defaultParamDim
$
such that 
\begin{enumerate}[label=(\roman*)]
\item 
it holds that $ \vartheta $ is the \Adam\ process 
with loss $ \defaultStochLoss $, 
learning rates $ ( \gamma_n )_{ n \in \N } $, 
batch sizes $ ( J_n )_{ n \in \N } $, 
momentum decay factors $ ( \alpha_n )_{ n \in \N } $, 
second moment decay factors $ ( \beta_n )_{ n \in \N } $, 
regularizing factor $ \varepsilon \in (0,\infty) $, 
initial value $ \Theta_0 $, and data $ ( X_{ n, j } )_{ (n,j) \in \N^2 } $ 
and 
\item 
it holds for all $ n \in \N $ that 
\begin{equation}
  \Theta_n 
  = 
  \avgPa_n \Theta_{ n - 1 } + ( 1 - \avgPa_n ) \vartheta_n
  .
\end{equation}
\end{enumerate}
\end{definition}

\begin{myalgorithm}{\Adam\ with geometrically weighted averaging}{alg:adam_avg_g}
\begin{algorithmic}[1]
\Statex\textbf{Setting:} The mathematical objects introduced in \cref{def:averaged_adam_geometric}
\Statex\textbf{Input:} 
$ \numberofsteps \in \N $
\Statex\textbf{Output:} 
Geometrically weighted averaged \Adam\ process $ \Theta_{ \numberofsteps } \in \R^\defaultParamDim $ after $ \numberofsteps $ steps 

\vspace{-2mm} 
\noindent\hspace*{-8.1mm}\rule{\dimexpr\linewidth+9.24mm}{0.4pt}  
\vspace{-3mm} 

\State $ \vartheta \gets \Theta_0 $
\State $ \theta \gets \Theta_0 $
\State $ \momentprocess \gets 0 $
\State $ \Momentprocess \gets 0 $
\For{$ \nos \in \cu{ 1, 2, \dots, \numberofsteps } $}
\State $ g \gets ( J_{ \nos } )^{ - 1 } \sum_{ j = 1 }^{ J_n } \defaultStochGradient( \vartheta, X_{ \nos, j } ) $
\State $ \momentprocess \gets \alpha_{ \nos } \momentprocess + ( 1 - \alpha_{ \nos } ) g $
\State $ \Momentprocess \gets \beta_{ \nos } \momentprocess + ( 1 - \beta_{ \nos } ) g^{ \otimes 2 } $ 
\mycomment{Square $ g^{ \otimes 2 } $ is understood componentwise}
\State $ \hat{\momentprocess} \gets \momentprocess / ( 1 - \prod_{ k = 1 }^{ \nos } \alpha_k ) $
\State $ \hat{\Momentprocess} \gets \Momentprocess / ( 1 - \prod_{ k = 1 }^{ \nos } \beta_k ) $
\State $ \vartheta \gets \vartheta - \gamma_{ \nos } \hat{\momentprocess} / ( \hat{\Momentprocess}^{ \otimes (1/2) } + \varepsilon ) $ 
\mycomment{Root $ \Momentprocess^{ \otimes ( 1 / 2 ) } $ is understood componentwise}
\State $ \theta \gets \avgPa_{ \nos } \theta + ( 1 - \avgPa_{ \nos } ) \vartheta $ 
\mycomment{Update averaged iterate}
\EndFor
\noindent 
\State \Return $ \theta $
\end{algorithmic}
\end{myalgorithm}

In the scientific literature the averaging procedure described in \cref{def:averaged_adam_geometric} and \cref{alg:adam_avg_g} above 
is typically referred to as \EMA. 
In particular, 
we refer, \eg, to \cite{morales-brotons2024exponential,SandlerZhmoginov2023} for numerical experiments for \EMA, 
we refer, \eg, to \cite{AhnCutkosky2024arXiv} for convergence analyses of \EMA\ in conjunction with \Adam, 
and we refer, \eg, to \cite{BusbridgeRamapuram2023EMA} for \SDE\ limits of \SGD\ methods with \EMA.

\section{Numerical experiments}
\label{sec:numerics}

\subsection{Introduction}
	We tested the two averaged variants of the \Adam\ optimizer in \cref{alg:adam_avg_a2} (see \cref{def:averaged_adam_arithmetic}) 
	and \cref{alg:adam_avg_g} (see \cref{def:averaged_adam_geometric}).
	Specifically, we employed \cref{alg:adam_avg_a2} with $ \AvgParam \in \cu{100 , 1000 } $ and 
	\cref{alg:adam_avg_g} with 
	$\forall \, n \in \N \colon \delta_n = \delta_1  \in \cu{0.99, 0.999}$,
	where $ \delta_1 = 0.999 $ seemed to work best in most cases.
	
	All of the experiments were implemented in the machine learning library \textsc{PyTorch} (cf., \eg, \cite{Paszke2017,Paszkeetal2019}).
	The {\sc Python} source codes for each of the performed numerical simulations 
	can be found on {\sc GitHub} at \href{https://github.com/deeplearningmethods/averaged-adam}{https://github.com/deeplearningmethods/averaged-adam}.

	\subsection{Polynomial regression}
	\label{subsec:poly_reg}
	
\begin{figure}
	\begin{center}  
		\includegraphics[scale=0.6]{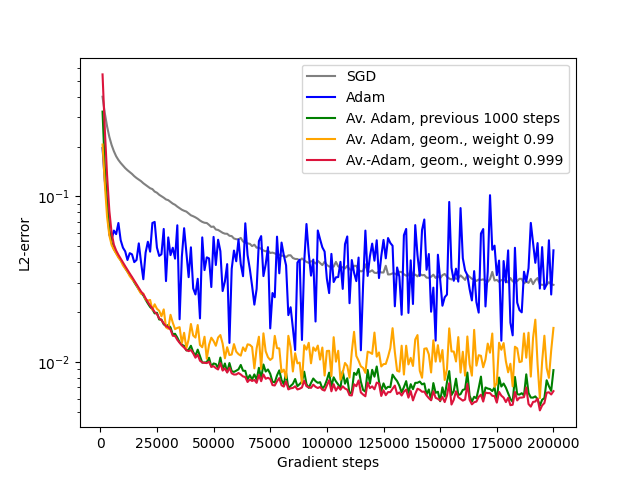}
		\caption{Numerical results for the polynomial regression problem described in \cref{subsec:poly_reg}.}
		\label{fig:polynomial}
	\end{center}
	\end{figure}

	As a first example we consider a simple regression problem with data corrupted by random noise.
	We optimize the coefficients of a polynomial with degree at most $ 25 $ to approximate the function
	$[-1, 1] \ni x \mapsto \sin ( \pi x ) \in \R$
	in $ L^2( [-1, 1]; \R ) $.
	In other words, we attempt to minimize the function 
	\begin{equation}
		\R^{d+1} \ni 
		\theta = ( \theta_0, \theta_1, \dots, \theta_d ) 
		\mapsto 
		\int_{ - 1 }^1 
		| \sin( \pi x ) - \ssum_{ k = 0 }^d \theta_k x^k |^2 \, \d x  
		\in \R
	\end{equation}
	for $d=25$, leading to an $26$-dimensional convex optimization problem.
	For the training we use a batch size of $64$ and constant learning rates of size $ 10^{ - 2 } $
	and we add random noise to the output following a centered Gaussian distribution with variance $ \nicefrac{ 1 }{ 5 } $. 
    Here and in most of the following numerical experiments we compare the plain vanilla \SGD\ method,
    the standard \Adam\ optimizer,
    \Adam\ with partial arithmetic averaging (\cref{alg:adam_avg_a2}) with $ \AvgParam = 1000 $, and 
    \Adam\ with geometrically weighted averaging (\cref{alg:adam_avg_g}) with 
    $\forall \, n \in \N \colon \delta_n = \delta_1  \in \cu{0.99, 0.999}$.
	The results are visualized in \cref{fig:polynomial}.
	
	\subsection{Artificial neural network (ANN) approximations for explicitly given functions}
	\label{subsec:simple_supervised}
	
	\subsubsection{ANN approximations for a 6-dimensional polynomial}
	We train standard fully connected feedforward \ANNs\ to approximate the target function
	\begin{equation}
		\label{eq:cubic_fct}
		\begin{split}
			[-1, 1 ] ^d \ni 
			x = ( x_1, \dots, x_d ) \mapsto 1 + \ssum_{i=1}^d (d + 1 - 2 i ) ( x_i )^3 \in \R
		\end{split}
	\end{equation}
	for $d= 6$.
	We use \ANNs\ with the \ReLU\ activation and two hidden layers consisting of $64$ neurons each.
	As the input distribution we choose the continuous uniform distribution on $ [ -1, 1]^d $.
	We use a batch size of 256 and constant learning rates of size $10^{-2}$.
	The results are visualized in \cref{fig:supervised}.
	
	\subsubsection{ANN approximations for a 20-dimensional normal distribution}
\begin{figure}
		\centering
		\begin{subfigure}{.5\textwidth}
			\centering
			\includegraphics[scale=0.55]{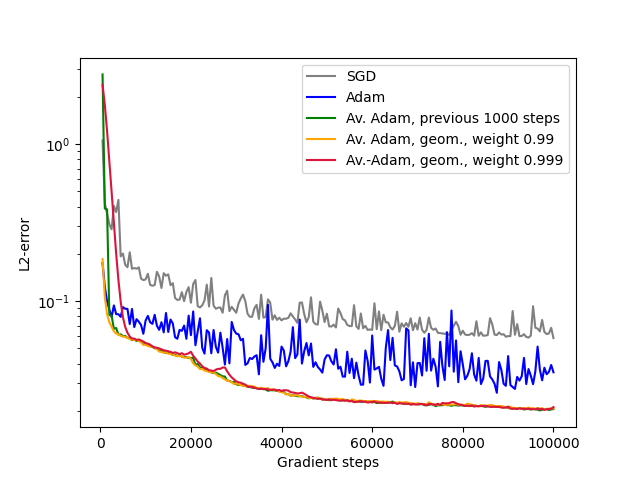}
		\end{subfigure}%
		\begin{subfigure}{.5\textwidth}
			\centering
			\includegraphics[scale=0.55]{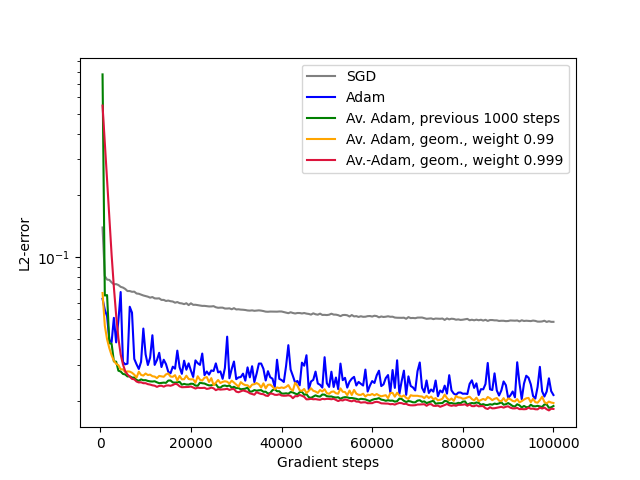}
		\end{subfigure}		
		\caption{Results for the supervised learning problem for the target functions in \cref{eq:cubic_fct} (left) and \cref{eq:gauss} (right).}
		\label{fig:supervised}
	\end{figure}

	As a further example we train fully connected feedforward \ANNs\ to approximate 
	the unnormalized density function 
	\begin{equation}
		\label{eq:gauss}
		[-2, 2]^d \ni x \mapsto 
		\exp\bigl( 
		  - \tfrac{ \norm{x}^2 }{ 6 } 
		\bigr) \in \R
	\end{equation}
	in $d = 20 $ dimensions.
	We use \ANNs\ with the \ReLU\ activation and three hidden layers consisting of $50$, $100$, and $50$ neurons, respectively.
	As the input distribution we choose the continuous uniform distribution on $ [-2, 2]^d $.
	For the training we employ the \Adam\ optimizer with a batch size of $256$ and constant learning rates of size $10^{-3}$.
	Again we add random noise to the output following a centered Gaussian distribution with variance $ \nicefrac{ 1 }{ 5 } $. 
	The results are visualized in \cref{fig:supervised}.

	\subsection{Deep Kolmogorov method (DKM)}
	\label{ssub:DKM}
	Within this subsection we use the \DKM\ proposed in Beck et al.~\cite{Beck2021Kolmogorov} to solve different linear \PDEs.
	\subsubsection{Heat PDE}
\begin{figure}[h]
		\centering 
		\begin{subfigure}{.49\textwidth}
			\centering 
			\includegraphics[scale=0.55]{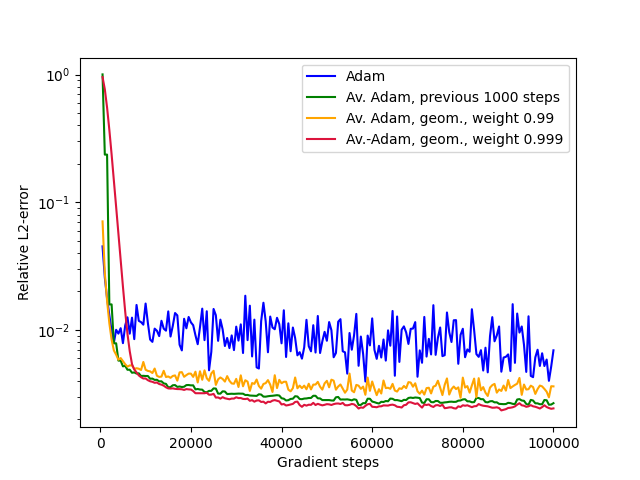}
		\end{subfigure}
		\begin{subfigure}{.49\textwidth}
			\centering 
			\includegraphics[scale=0.55]{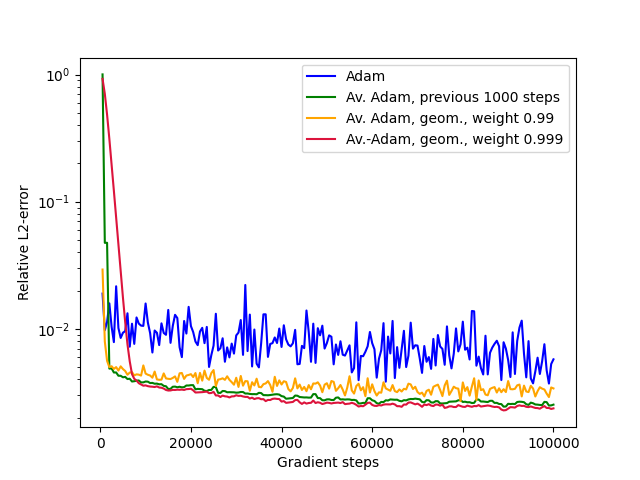}
		\end{subfigure}
		\caption{Results for the heat \PDE\ in \cref{eq:heat} using the \DKM with constant learning rates (left) and decreasing learning rates (right).}
		\label{fig:heat}
	\end{figure}

	We first consider the heat \PDE\ on $\R ^d$ for 
	$d=10$.
	Specifically, we consider $ d = 10 $ and we attempt to approximate the solution $u \colon [0 , T] \times \R^d \to \R$ of the \PDE
	\begin{equation}
		\label{eq:heat}
			\tfrac{ \partial u }{ \partial t } = \Delta_x u ,
		\qquad 
			u(0, x ) 
			= 
			\norm{x}^2
	\end{equation}
	for $ t \in [0,T] $, $ x \in \R^d $ 
	at the final time $T=2$ on the domain $[-1, 1] ^d$.
	The \PDE\ can be reformulated as a stochastic minimization problem (cf.~Beck et al.~\cite{Beck2021Kolmogorov}) 
	and thus \SGD\ methods such as the \Adam\ optimizer can be used to compute an approximate minimizer.
	We employ fully connected feedforward \ANNs\ with three hidden layers consisting of $50$, $100$, and $50$ neurons, respectively. 
	This time we use the smooth \GELU\ activation, which seems to be more suitable for \PDE\ problems.
	We tested both constant learning rates of size $5 \cdot 10^{-4}$
	and polynomially decaying learning rates of size 
	$\gamma_n = 5 \cdot 10^{-3} \cdot  n^{ - 1 / 4}$, and a batch size of $2048$.
	To compute the test error we compare the output with the exact solution
	$ u(t,x) = \norm{x}^2 + 2 d t $ 
	for $ t \in [0,T] $, $ x \in \R^d $.
	The results are visualized in \cref{fig:heat}.
	
	\subsubsection{Black--Scholes PDE}

\begin{figure}
		\centering
	\begin{subfigure}{0.49\textwidth}
		\centering 
		\includegraphics[scale=0.55]{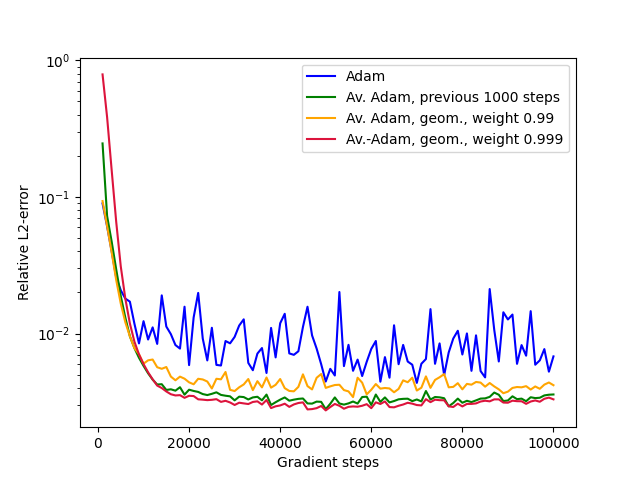}
	\end{subfigure}
	\begin{subfigure}{0.49\textwidth}
		\centering
		\includegraphics[scale=0.55]{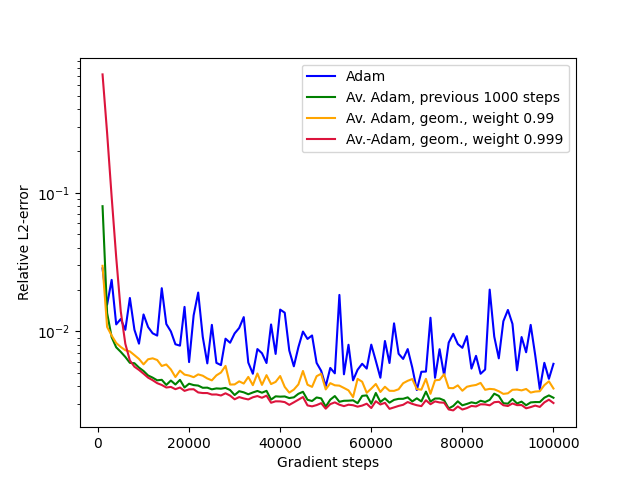}
	\end{subfigure}
\caption{Results for the Black-Scholes \PDE\ in \cref{eq:black-scholes} using the \DKM,
	with constant learning rates (left) and decreasing learning rates (right).}
	\label{fig:bs}
	\end{figure}

	We next consider a Black--Scholes \PDE\ on $ [90, 110]^d $ for $ d = 20 $.
	Specifically, we consider $ d = 20 $ and 
	we attempt to approximate the solution $ u \colon [0, T] \times \R^d \to \R $ of the \PDE
	\begin{equation}
		\label{eq:black-scholes}
\textstyle 
			\tfrac{\partial u }{\partial t }
			= \frac{1}{2} \sum\limits_{i=1}^d \abs{ \sigma_i x_i }^2 \tfrac{\partial^2 u }{ \partial x_i^2 } 
			+ \mu \sum\limits_{ i = 1 }^d x_i \tfrac{ \partial u }{ \partial x_i } ,
			\quad 
			u( 0 , x ) 
			= \exp( - r T ) \max \cu*{ \max \cu{ x_1, x_2, \dots, x_d } - K , 0 }
	\end{equation}
	for $ t \in [0,T] $, $ x = ( x_1, \dots, x_d ) \in \R^d $ 
	where $ \sigma = ( \sigma_i )_{ i \in \{ 1, 2, \dots, d \} } = ( \frac{ i + 1 }{ 2 d } )_{ i \in \{ 1, 2, \dots, d \} } $, 
	$ r = - \mu = \frac{ 1 }{ 20 } $, $ K = 100 $ 
	at the final time $T=1$ using the \DKM. 
	We again employ fully connected feedforward \ANNs\ with the \GELU\ activation and three hidden layers consisting of $50$, $100$, and $50$ neurons, respectively,
	and a batch normalization layer before the first hidden layer.
	For the training we use the batch size $2048$ and test two different learning rate schedules:
	Constant learning rates of size $5 \cdot 10^{-4}$ and slowly decreasing learning rates of the form $\gamma_n = 5 \cdot 10^{-3} \cdot n^{- \nicefrac{1}{4} }$ (see \cref{fig:bs}), which lead to comparable results.
	To compute the test error we compare the output 
	with the exact solution computed with the Feynman--Kac formula and approximated using a Monte Carlo method with $2048000$ Monte Carlo samples.

\begin{figure}
	\begin{center}   
		\includegraphics[scale=0.55]{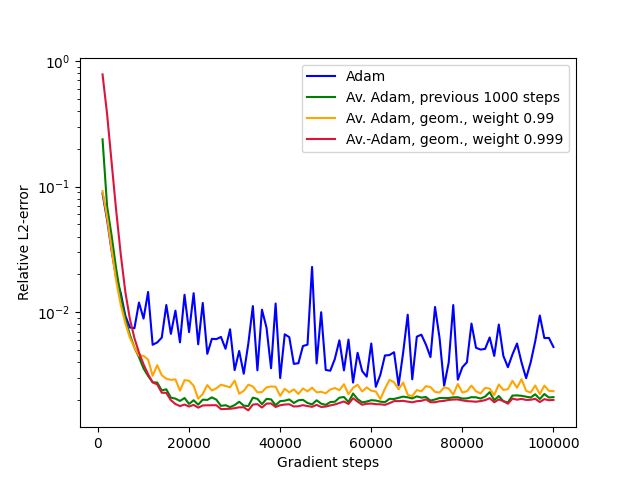}
		\caption{Results for the Black-Scholes \PDE\ in \cref{eq:black-scholes-cor} using the \DKM.}
		\label{fig:bs_cor}
	\end{center}
	\end{figure}

	As a second example we consider a Black--Scholes \PDE\ with correlated noise on $[90, 110] ^d$ for $d=20$.
	Specifically, let $ d = 20 $ and 
	let $ Q = ( Q_{ i, j } )_{ (i,j) \in \{ 1, 2, \dots, d \}^2 } $, 
	$ \Sigma = ( \Sigma_{ i, j } )_{ (i,j) \in \{ 1, 2, \dots, d \}^2 } \in \R^{ d \times d }$, 
	$ \beta = ( \beta_1, \dots, \beta_d ) $, $ \zeta_1, \zeta_2, \dots, \zeta_d \in \R^d$
	satisfy for all
	$i, j , k \in \cu{1, 2, \ldots, d }$ with $i < j$ that
	$\beta_k = \frac{1}{10} + \frac{k}{2 d }$,
	$Q_{k, k} = 1$,
	$Q_{i, j } = Q_{j, i } = \frac{1}{2}$,
	$\Sigma_{i , j} = 0$,
	$\Sigma_{k, k } > 0$,
	$\Sigma \Sigma^* = Q$,
	and
	$
	  \zeta_k = ( \beta_k \Sigma_{ k, 1 }, \beta_k \Sigma_{ k, 2 }, \dots, \beta_k \Sigma_{ k, d } )
	$
	(cf.~\cite[Section 4.4]{Beck2021Kolmogorov}). 
	We attempt to approximate the solution $u \colon [0, T] \times \R^d \to \R$ of the \PDE
	\begin{equation}
		\label{eq:black-scholes-cor}
		\begin{split}
		\textstyle 
			\tfrac{\partial u }{\partial t }
			&
		\textstyle 
			= \frac{1}{2} \sum\limits_{ i, j = 1 }^d x_i x_j \spro{\zeta_i, \zeta_j} \tfrac{\partial ^2 u }{ \partial x_i \partial x_j } 
			+ \mu \sum\limits_{ i = 1 }^d x_i \tfrac{\partial u }{\partial x_i} , \\
			u ( 0 , x ) 
			&= \exp( - r T ) \max\cu*{ K - \min \cu{ x_1, x_2, \dots, x_d } , 0 } 
		\end{split}
	\end{equation}
	for $ t \in [0,T] $, $ x = ( x_1, \dots, x_d ) \in \R^d $
	where $r = - \mu = \frac{1}{20}$, $K=110$ at the final time $T=1$ using the \DKM. 
	The remaining hyperparameters for the experiment are the same as for the Black--Scholes \PDE\ in \cref{eq:black-scholes}.
	The results are visualized in \cref{fig:bs_cor}.
	
	\subsection{Physics-informed neural networks (PINNs)}
	\label{ssub:PINNs}
	We used the method of \PINNs\ to approximately solve a few semilinear heat \PDEs.
	
	\subsubsection{Allen-Cahn PDE}
\begin{figure}
		\centering
	\begin{subfigure}{.49\textwidth}
		\centering 
		\includegraphics[scale=0.55]{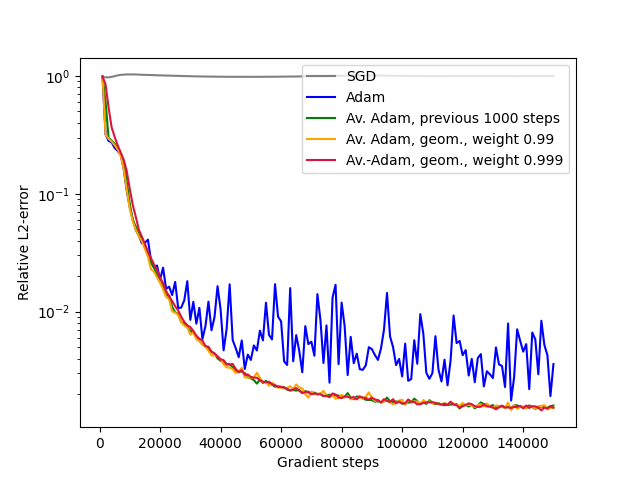}		
	\end{subfigure}
	\begin{subfigure}{.49\textwidth}
		\centering
		\includegraphics[scale=0.55]{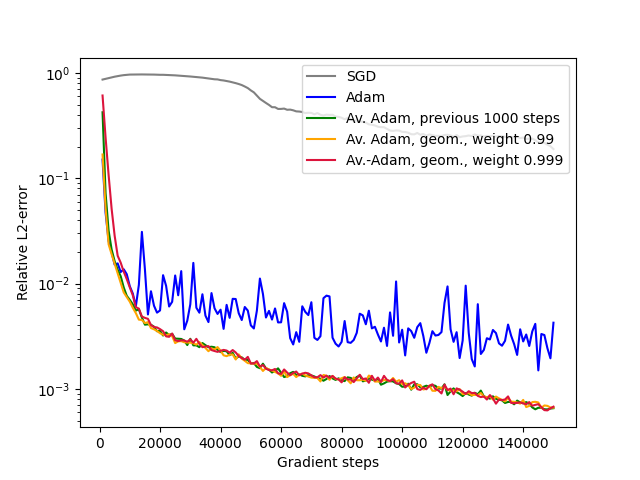}
	\end{subfigure}
	\caption{Results for the Allen-Cahn \PDE\ in \cref{eq:pinn-ac} (left)
		and the Sine-Gordon \PDE\ in \cref{eq:pinn-sine} (right) using \PINNs.}
	\label{fig:pinn}
	\end{figure}

	We first consider the Allen-Cahn \PDE\ on the domain $D = [0,2] \times [0,1]$ with time horizon $T=4$. 
	Specifically, we attempt to approximate the solution $ u \colon [0, T] \times D \to \R$ of the \PDE\
	\begin{equation}
		\label{eq:pinn-ac}
		\begin{split}
			\tfrac{ \partial u }{ \partial t } = \tfrac{1}{100} \Delta_x u + (u - u^3) ,
			\qquad 
			u(0, x ) = 
			\sin(\pi x_1) \sin(\pi x_2)
		\end{split}
	\end{equation}
	for $ t \in [0,T] $, $ x = ( x_1, x_2 ) \in D $
	equipped with Dirichlet boundary conditions at the terminal time $ T = 4 $.
	We used fully connected feedforward \ANNs\ with 
	3 hidden layers consisting of $32$, $64$, and $32$ neurons, respectively, and the \GELU\ activation.
	For the training we employ the \Adam\ optimizer with a batch size of $256$ and constant learning rates of size $10^{-3}$.
	To compute the test error we compare the output with a reference solution obtained by a finite element method using $101^2$ degrees of freedom in space and $500$ second order linear
	implicit Runge-Kutta time steps.
	The results are visualized in \cref{fig:pinn}.
	
	\subsubsection{Sine-Gordon type PDE}

	We next consider a Sine-Gordon type semilinear \PDE\ on the domain $D = [0,2] \times [0,1]$ with time horizon $T=1$.
	Specifically, 
	we attempt to approximate the solution $ u \colon [0, T] \times D \to \R$ of the \PDE\
	\begin{equation}
	\textstyle 
		\label{eq:pinn-sine}
			\frac{ \partial u }{ \partial t } = \tfrac{1}{20} \Delta_x u + \sin (u) , 
			\qquad 
			u(0, x ) = \tfrac{3}{2} | \sin( \pi x_1 ) \sin( \pi x_2 ) |^2
	\end{equation}
	for $ t \in [0,T] $, $ x = ( x_1, x_2 ) \in D $
	equipped with Dirichlet boundary conditions at the terminal time $ T = 1 $.
	The other hyperparameters for the training are the same as in the case of the Allen-Cahn \PDE.
	The results are visualized in \cref{fig:pinn}.

	\subsubsection{Burgers equation}
	
	We also employ the \PINN\ method to approximately solve the one-dimensional Burgers equation
	\begin{equation}
	\textstyle 
		\label{eq:pinn-burgers}
			\frac{ \partial u }{ \partial t } = \alpha \Delta_x u - u \frac{ \partial u }{ \partial x } 
			,
	\qquad 
			u(0, x ) = \frac{ 2 \alpha \pi \sin( \pi x ) }{ \beta + \cos( \pi x ) } 
	\end{equation}
	for $ t \in [0,T] $, $ x \in D = [0,2] $
	equipped with Dirichlet boundary conditions at the terminal time $T = \frac{1}{2}$, 
	where $\alpha = \frac{1}{20}$ and $\beta = \frac{11}{10}$.
	The exact solution satisfies for all $ t \in [0,T] $, $ x \in D $ that
	$ u(t,x) = \frac{ 2 \alpha  \pi \sin ( \pi x ) }{ \beta \exp ( \alpha t \pi ^2 ) + \cos ( \pi x ) } $.
	
\begin{figure}
	\begin{center}
		\includegraphics[scale=0.5]{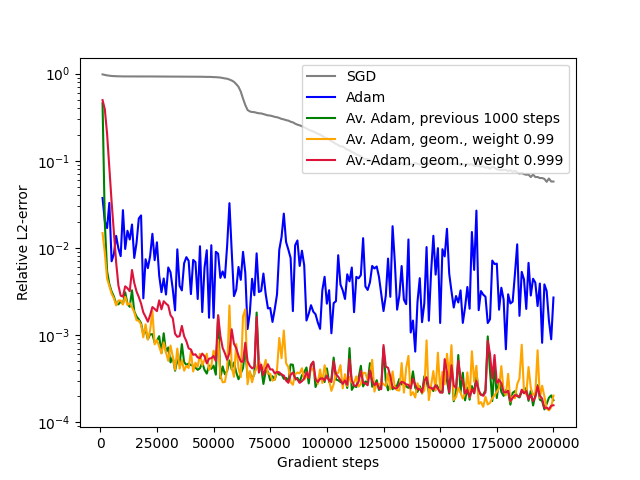}
		\caption{Results for the Burgers \PDE\ in \cref{eq:pinn-burgers} using \PINNs.}
		\label{fig:pinn_burg}
	\end{center}
	\end{figure}

	We use fully connected feedforward \ANNs\ with 
	3 hidden layers consisting of $16$, $32$, and $16$ neurons, respectively, and the \GELU\ activation.
	For the training we employ the \Adam\ optimizer with a batch size of $128$ and constant learning rates of size $3 \cdot 10^{-3}$.
	The results are visualized in \cref{fig:pinn_burg}.
	
	\subsection{Deep BSDE method for a Hamiltonian--Jacobi--Bellman (HJB) equation}
	\label{ssub:deepBSDE}
	We employ the \deepBSDE\ method introduced in E et al.~\cite{EHanJentzen2017,HanJentzenE2018}
	to approximate the solution $u \colon [0, T ] \times \R^d \to \R$
	of the Hamilton--Jacobi--Bellman \PDE 
	\begin{equation}
		\label{eq:bsde}
			\tfrac{ \partial u }{ \partial t } = - \Delta_x u + \norm{ \nabla_x u }^2 ,
			\qquad 
			u(T, x ) = \ln \rbr*{ \tfrac{1}{2} ( \norm{x}^2 + 1 ) }
	\end{equation}
for $ t \in [0,T] $, $ x \in \R^d $ at initial time $0$ for the time horizon $ T = \nicefrac{ 1 }{ 4 } $
and dimension $d = 25$.
We approximately solve the \PDE\ on the domain $[-1, 1] ^d$ using an \ANN\ with 
two hidden layers consisting of 45 neurons each and the \GELU\ activation.
We employ a time discretization with $N=20$ time steps
and approximate the gradient of the solution at each time step using an \ANN\ with 
two hidden layers consisting of 45 neurons each and the \GELU\ activation.

\begin{figure}
	\begin{center}
		\includegraphics[scale=0.5]{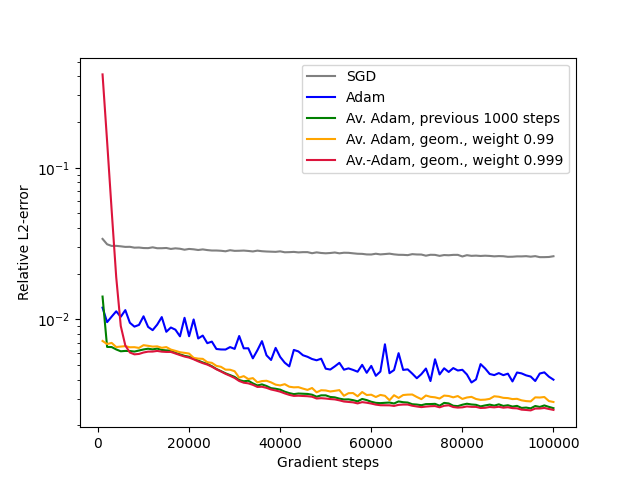}
		\caption{Results for the \PDE\ in \cref{eq:bsde} using the \deepBSDE\ method.}
		\label{fig:deep-bsde}
	\end{center}
	\end{figure}

For the training we employ the \Adam\ optimizer with a batch size of 512 and 
slowly decreasing learning rates of the form $\gamma_n = 0.02 \cdot  n^{ - 1 / 5 }$.
To compute the test error we compare the output with an approximation of the exact solution 
computed through a Monte Carlo approximation with $ 819200 $ Monte Carlo samples
for the Cole--Hopf transform (cf., \eg, \cite[Lemma 4.2]{EHanJentzen2017}). 
The results are visualized in \cref{fig:deep-bsde}.

	\subsection{Optimal control problem}
	\label{ssub:OC}
	We consider a controlled diffusion of the form
	\begin{equation}
		\label{eq:control_sde}
		\d X_t = (A X_t + B u_t ) \, \d t + \sqrt{2} \, \d W_t
	\end{equation}
	for $ t \in [0,\infty) $ where $ d \in \N $, 
	where $ X \colon [0, \infty ) \times \Omega \to \R^d$ is the diffusion process,
	where $u \colon [0, \infty ) \times \Omega \to \R^d$ is the control process,
	and where $A, B \in \R^{d \times d }$ are given matrices.
	We define the cost functional
	\begin{equation}
		J( t , x, u ) 
		=
		\E\br*{ 
		  \int_t^T \rbr*{\tfrac{1}{2} \norm{X_s} ^2 + \norm{u_s } ^2 } \, \d s + \norm{X_T } ^2 \Big| X_t = x 
		}
	\end{equation}
	and attempt to compute the minimal expected cost
	\begin{equation}
		\label{eq:control_inf}
		\inf_u \E\bigl[ J ( 0 , Z , u ) \bigr] ,
	\end{equation}
	where $Z \colon \Omega \to [-1, 1 ] ^d$ is assumed to be continuous uniformly distributed. 
	To this end, let $ \timesteps \in \N $, 
	$ t_0, t_1, \dots, t_{ \timesteps } \in \R $ satisfy 
	for all $ k \in \{ 1, 2, \dots, \timesteps \} $ that 
	$ t_k = \frac{ k T }{ N } $.
	We approximate the solution of the \SDE\ in \cref{eq:control_sde} with a forward Euler method. 
	For each $n \in \cu{0, 1, \ldots, N-1}$ we approximate the control $ u_{ t_n } $ through $ u_{ t_n } \approx \realization{\theta_n}( X_{ t_n } ) $ 
	where $ \realization{\theta_n} \colon \R^d \to \R^d$ is the realization function of a neural network with parameter vector $\theta_n$.
	In our experiment we use the values $d=10$, $T=1$, $N=100$, and \ANNs\ with 2 hidden layers consisting of 40 neurons each and the \GELU\ activation.
	Additionally, we employ batch normalization after the input layer and each hidden layer.
	We train these neural networks using the \Adam\ optimizer with a batch size of $1024$ and constant learning rates of size $10 ^{ - 2}$.

	To compute the test error we approximate the expectation in \cref{eq:control_inf} with $4096$ Monte Carlo samples for $Z$ and $500$ Monte Carlo samples for the Brownian motion $W$
	and compare the output with the exact solution obtained by solving the corresponding Ricatti \ODE\ with $100000$ time steps.
	
\begin{figure}
	\begin{center}
		\includegraphics[scale=0.5]{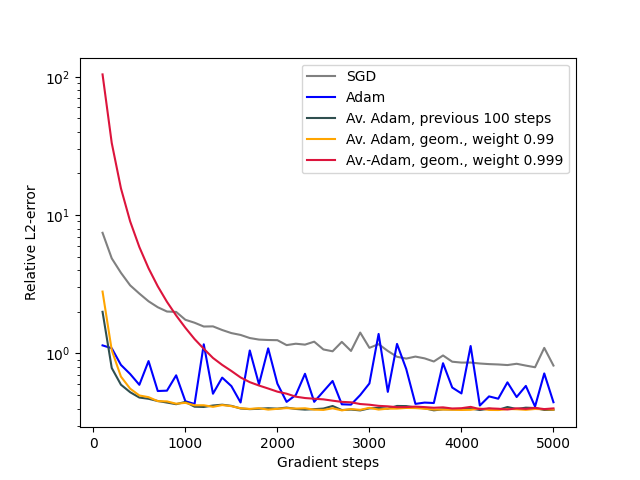}
		\caption{Results for the optimal control problem in \cref{eq:control_inf}.}
		\label{fig:lqr}
	\end{center}
	\end{figure}

	This time we compare the plain vanilla \SGD\ method,
	the standard \Adam\ optimizer,
	\Adam\ with partial arithmetic averaging (\cref{alg:adam_avg_a2}) with $ \AvgParam = 100 $, and 
	\Adam\ with geometrically weighted averaging (\cref{alg:adam_avg_g}) with  
	$ \forall \, n \in \N \colon \delta_n = \delta_1 \in \cu{ 0.99 , 0.999 }$.
	The results are visualized in \cref{fig:lqr}.

	\subsection{Image classification}
	\label{ssub:CIFAR-10}
\begin{figure}
	\begin{center}   
		\includegraphics[scale=0.5]{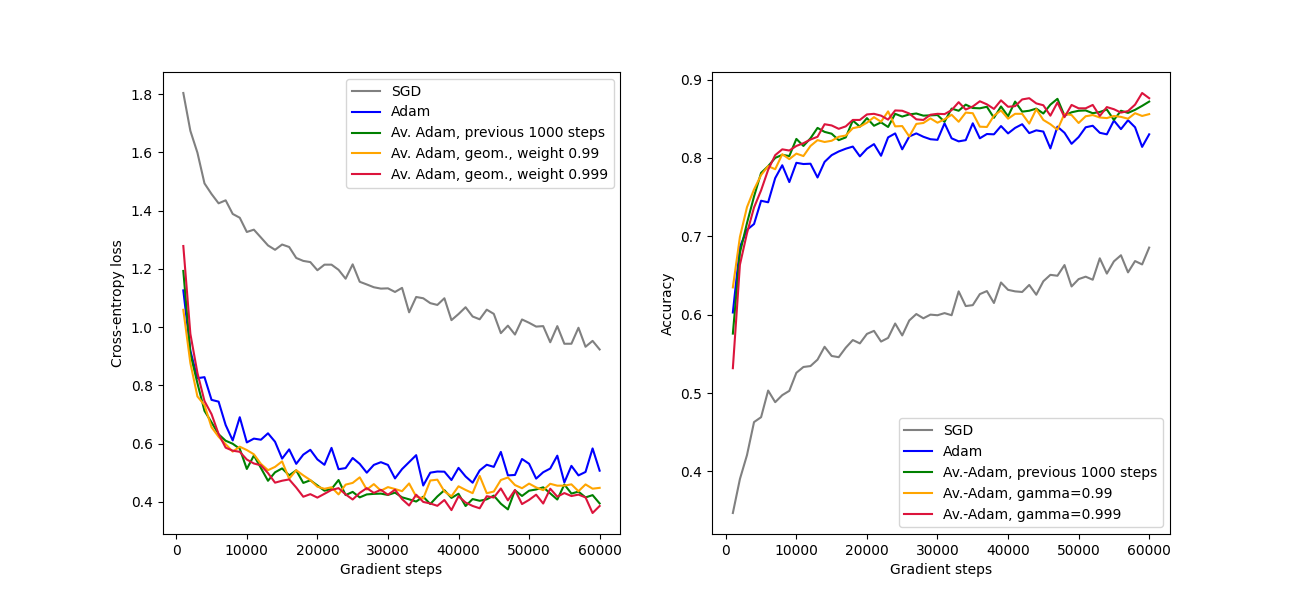}
		\caption{Results for the CIFAR-10 image classification problem using a \ResNet.
		Left: test error, right: test accuracy.}
		\label{fig:cifar}
	\end{center}
	\end{figure}

	We train a \ResNet\ on the CIFAR-10 dataset, a standard benchmark problem for image classification.
	Specifically, we use a variant of the \ResNet\ architecture described in He et al.~\cite{He_2016_CVPR} with 9 layers, following the implementation at \url{https://github.com/matthias-wright/cifar10-resnet}.
	We employ the \Adam\ optimizer with a batch size of $64$,
	constant learning rates of size $2 \cdot 10 ^{ - 4 }$, and weight decay with decay parameter $3 \cdot 10^{-4}$.
	Following standard ideas from data augmentation,
	we apply random horizontal flips with probability $ \nicefrac{1}{2} $,
	random offset cropping down to $32 \times 32$, using reflection padding of 4 pixels,
	random color jitter,
	and random rotations by an angle sampled uniformly from the interval $ ( 0, \nicefrac{ \pi }{ 4 } ) $, to the training data.
	For details we refer to the documentation of {\sc PyTorch} transforms.
	The results are visualized in \cref{fig:cifar}.
	
\section{Conclusion}
\label{sec:conclusion}

In this work we have applied different averaged variants of the \Adam\ optimizer 
(see \cref{sec:methods}) to a series of learning problems 
including polynomial regression problems 
(see \cref{subsec:poly_reg}), 
including deep \ANN\ approximations 
for explicitly given high-dimensional target functions
(see \cref{subsec:simple_supervised}), 
including deep \ANN\ approximations for stochastic \OC\ problems 
(see \cref{ssub:OC}), 
including 
\DK\ (see \cref{ssub:DKM}), 
\PINN\ (see \cref{ssub:PINNs}), and 
\deepBSDE[\emph{deep backward stochastic differential equation}]\ (see \cref{ssub:deepBSDE}) 
approximations 
for \PDEs, 
and 
including residual deep \ANN\ approximations 
for the CIFAR-10 image classification dataset 
(see \cref{ssub:CIFAR-10}). 
%
%
%
%
In each of the considered numerical examples 
the employed averaged variants of the \Adam\ optimizer 
outperform the standard \Adam\ and the standard \SGD\ optimizers, 
particularly, in the situation of 
the considered scientific computing problems.
Taking this into account, 
we believe that it is very relevant to further study and employ averaged variants 
of the \Adam\ and similar optimizers, particularly, 
when solving \PDE, \OC, or related scientific computing problems by 
means of deep learning approximation methods.

\subsubsection*{Acknowledgments}

This work has been partially funded by the European Union (ERC, MONTECARLO, 101045811). 
The views and the opinions expressed in this work are however those of the authors only and do 
not necessarily reflect those of the European Union or the European Research Council (ERC). 
Neither the European Union nor the granting authority can be held responsible for them. 
In addition, this work has been partially funded by the Deutsche Forschungsgemeinschaft (DFG, German Research
Foundation) under Germany’s Excellence Strategy EXC 2044-390685587, Mathematics Münster: Dynamics-Geometry-Structure.
Moreover, this work has been partially supported 
by the Ministry of Culture and Science NRW 
as part of the Lamarr Fellow Network.


\end{document}